\documentclass[preprint,amsmath,amssymb,
 aps]{revtex4-1}

\usepackage{graphicx}
\usepackage{dcolumn}
\usepackage{bm}

\DeclareSymbolFont{AMSb}{U}{msb}{m}{n}
\DeclareSymbolFontAlphabet{\Bbb}{AMSb}
\DeclareMathSymbol{\bba}{\mathbin}{AMSb}{'101}
\DeclareMathSymbol{\bbr}{\mathbin}{AMSb}{'122}
\DeclareMathSymbol{\bbc}{\mathbin}{AMSb}{'103}
\DeclareMathSymbol{\bbh}{\mathbin}{AMSb}{'110}
\DeclareMathSymbol{\bbo}{\mathbin}{AMSb}{'117}

\begin{document}

\preprint{kissing2013}

\title{
Delimiting Maximal Kissing Configurations in Four Dimensions
}

\author{Eric Lewin Altschuler}
\affiliation{
Departments of Physical Medicine and Rehabilitation and Microbiology \& Molecular Genetics,
New Jersey Medical School\\ University Hospital, B-403, 150 Bergen Street,
 Newark, NJ 07103, USA}
 \email{ eric.altschuler@umdnj.edu}

\author{Antonio P\'erez--Garrido}
\affiliation{
Departamento de F\'\i sica Aplicada, UPCT
Campus Muralla del Mar, Cartagena, 30202 Murcia, Spain}
\email{antonio.perez@upct.es}
 
\date{\today}

\begin{abstract}
How many unit $n-$dimensional spheres can simultaneously touch or kiss a central $n-$dimensional unit sphere? Beyond mathematics this question has implications for fields such as cryptography and the structure of biologic and chemical macromolecules. The kissing number is only known for dimensions 1-4, 8 and 24  (2, 6, 12, 24, 240, 19650, respectively) and only particularly obvious for dimensions one and two. Indeed, in four dimensions it is not even known if Platonic polytope unique to that dimension known as the 24-cell is the unique kissing configuration.
We have not been able to prove that the 24-cell is unique, but, using a physical approach utilizing the hopf map from four to three dimensions, we for the first time delimit the possible other configurations which could be kissing in four dimensions.
\end{abstract}
\pacs{PACS}

\maketitle

How many solid unit $n$-dimensional spheres can be placed such that each touches a central unit $n$-dimensional sphere?  This question is known as the
$n$-dimensional kissing problem, and the number of unit spheres that can so touch the central sphere the $n$-dimensional kissing number.  
The kissing problem has applications to many and diverse fields: optical systems\cite{Me98}, cryptography\cite{CS98}, botany\cite{Ta30}, and understanding the structure of biologic 
\cite{BG03}
and chemical macromolecules\cite{LD03}.  As well, physics inspired/physical based thinking about arrangement problems\cite{MR08}
 can then lead to observation of new physical phenomena\cite{BM12}. 
In one-dimension the unit sphere is a closed unit line segment, and the one-dimensional kissing number is clearly two.  In two dimension the unit sphere is unit disc and the kissing number (as one can readily see by using seven of the same round coin)  is six.  In three dimensions the problem is more subtle.  Indeed, in a famous dispute (possibly apocryphal\cite{Ca11}) in the 1690's Isaac Newton argued that the kissing number is 12, 
while David Gregory thought it was 13. Perhaps not surprisingly, Newton was correct, though it took more than two and a half centuries to prove it\cite{SW53}.  In higher dimensions the problem becomes even harder. Applying linear programming methods of Delsarte \cite{De72,De73,DG77} to exceptionally symmetric structures that exist in 8 and 24 dimensions
Levenstein \cite{Le79} and Odlyzko and Sloane \cite{OS79} were able to prove in 1979 that the kissing numbers are 240 and 196560, respectively.  
Despite the existence of a plantonic polytope (known as the 24-cell) unique to  four dimensions which is a kissing configuration of 24 hyperspheres, proving that 24 is indeed the kissing number in four dimensions was a more stubborn problem than for 8 or 24 dimensions.  Indeed, the Delsarte method as typically applied was shown \cite{AB00} to bound the kissing number only to 25 or less.  A few years ago using ingenious and extensive applications of Delsarte's method Musin\cite{Mu03} and subsequently others using semidefinite programmings\cite{BV09,MV09} were able to prove that in fact the kissing number in four dimensions is 24.  But it is still not even known if the 24-cell is the only configuration in 4-dimensions with kissing number of 24.

\begin{figure}

\leavevmode
\begin{center}
\includegraphics[angle=0,width=8cm]{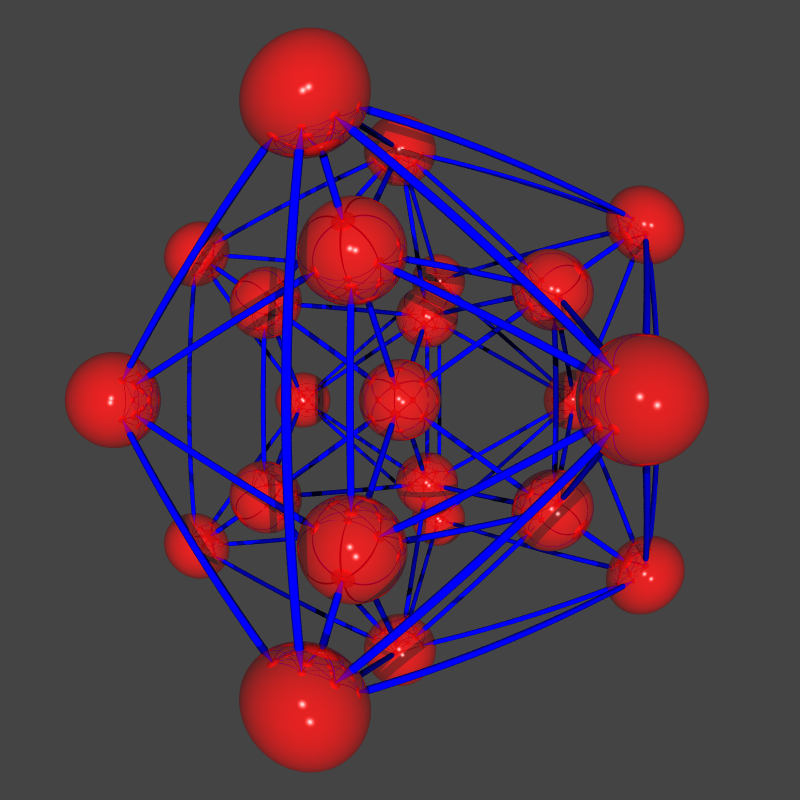}
\caption{24-cell in a Schlegel-like representation}
\label{24cell}
\end{center}
\end{figure}

The 24-cell is shown in Figure \ref{24cell}.  Our study of kissing configurations in 4-dimensions is aided by what is known as the 
Hopf map from four to three dimensions.   The Hopf map takes  points in four dimensions $(w, z)$ ---with the four coordinates 
as the components of the complex numbers $w$ and $z$---and $|w|^2+|z|^2=1$ on the surface of a four dimensional sphere ($S^3$) to a pair of numbers $(w, z)$ in 
$\bbc\times\bbr=\bbr^3$
\begin{equation}
(w,z) \rightarrow (2 w z^*, |z|^2-|w|^2)\,  {\rm in}\,\,  \bbc\times\bbr = \bbr^3.
\end{equation}
We check:
\begin{equation}
|2 w z^*|^2 + (|z|^2-|w|^2)^2
= 4 |w|^2 |z|^2 + (|z|^2-|w|^2)^2 = (|z|^2+|w|^2)^2 = 1.
\end{equation}
So this does map $S^3$ to the ordinary sphere $S^2$. If one fixes a point of the
ordinary sphere say  $(a,t)$  where  $a$  is complex, $t$ is real and
$|a|^2 + t^2 = 1$, then what its known as its {\it fiber}, i.e., the set of all points
which map to it, is a circle
\begin{equation}
\left(\frac{a e^{i \theta}}{\sqrt{2(1+t)}}, e^{i \theta}\sqrt{(1+t)/2}\right).
\label{map3-4}
\end{equation}
Further details, discussions and proof 
of the Hopf map from $S^3$ to $S^2$ is given in \cite{Ho31}.  

\begin{figure}
\leavevmode
\begin{center}
\includegraphics[angle=0,width=12cm]{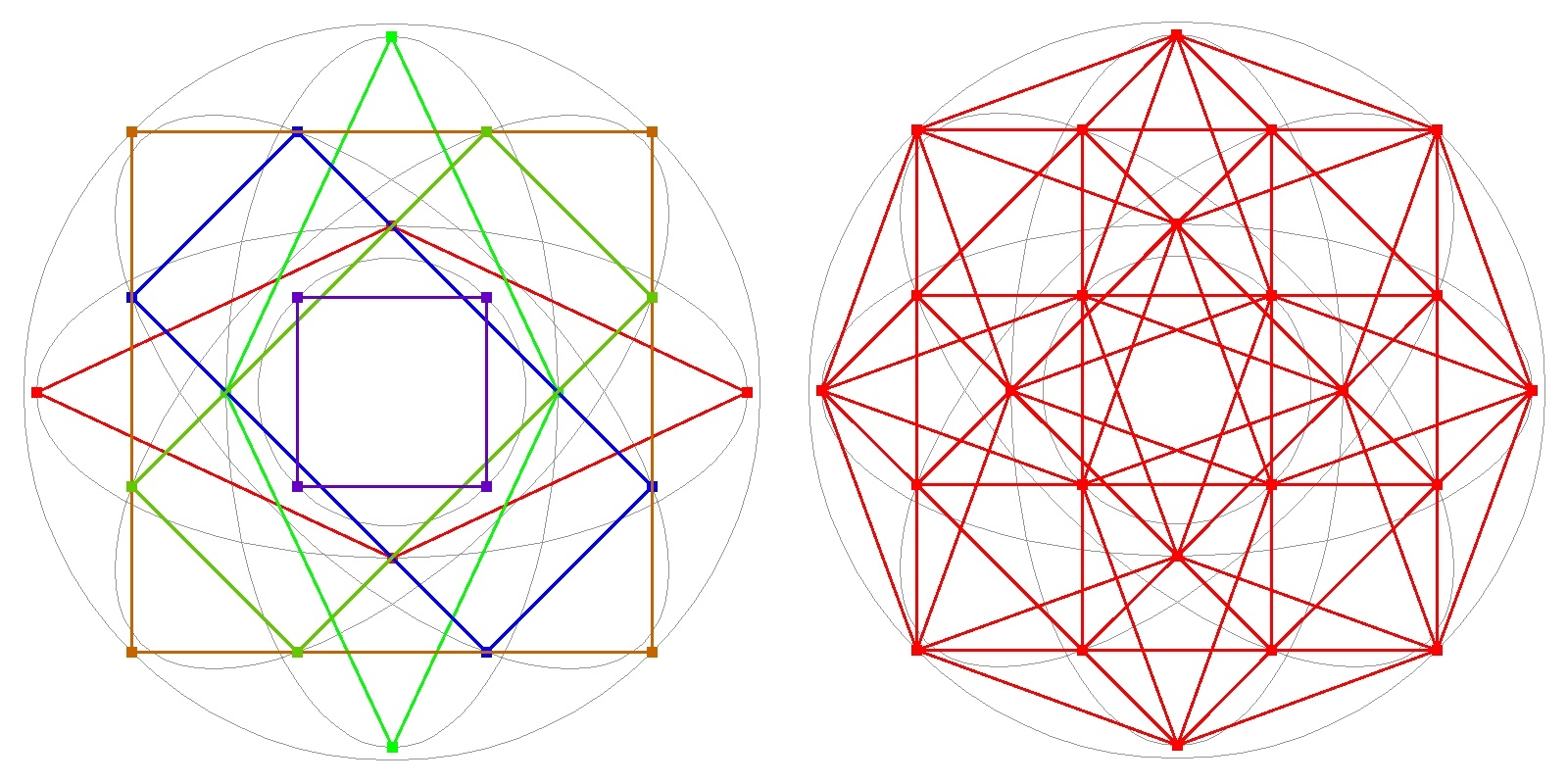}
\caption{24-cell in a orthographic projection to a plane. In the left figure lines join points on the same circle, each circle is drawn in a diferent color. Gray lines 
correspond to full circles. In the right figure lines join nearest neighbors points.}
\label{6x4}
\end{center}
\end{figure}

Intuitively or physically one can think of coordinates points on the surface of a unit sphere in four dimensions ($S^3$) as two spherical polar coordinates on a 2-sphere ($S^2$) (the surface of a 3-dimensional sphere) and the third coordinate being an azimuthal angle around a circle at the point on the 2-sphere. Since a point on 
$S^2$ lifts to a circle on $S^3$, henceforth points on $S^2$ will be called circles, e.g.,  when we say a circle on the north pole, we mean a point on the north pole 
of $S^2$ that after Hopf fibration becomes a circle on $S^3$. We call kissing points those on $S^3$ separated by a distance larger or equal to 1. A representation of the 
24-cell as six circles each with four points on it is shown in Figure \ref{6x4}. This ``Hopf perspective'' of the 3-sphere gives a simple appreciation of why the 24-cell is a kissing configuration\cite{AP07}.

We derive a relation to obtain distances on the 3-sphere ($d_3$) from polar coordinates on the 2-sphere and azimuthal angle on a circle ($\theta_i$). 
If we have two circles on $S^2$ separated by $d_2$ and lift them to $S^3$ using the Hopf map, we have that
distance on $S^3$ ($d_3$) is given by:
\begin{equation}
d^2_3=2-\sqrt{4-d^2_2}\cos\left( \theta_j-\theta_i+ \Phi_{ij}\right),
\label{d3}
\end{equation}
where $\Phi_{ij}$ is an angle that depends on original coordinates in $S^2$:
\begin{equation}
\cos \Phi_{ij}=\frac {2}{\sqrt{4-d_3^2}}\left(\cos (\phi_j-\phi_i)\sin \frac{\alpha_i}{2} \sin \frac{\alpha_j}{2}+ 
   \frac{\alpha_i}{2} \cos \frac{\alpha_j}{2}    \right),
\end{equation}
and
\begin{equation}
\sin   \Phi_{ij}=\frac {2}{\sqrt{4-d_3^2}} \sin(\phi_j-\phi_i)\sin \frac{\alpha_i}{2} \sin \frac{\alpha_j}{2},
\end{equation}
where $\phi_i$ and $\alpha_i$ are the polar coordinates of  point $i$ on $S^2$ ($\alpha\in \left[0-\pi\right]$, $\phi\in \left[0-2\pi\right]$).
We observe that when one circle is on the north pole then $\Phi_{ij}=0$ no matter where we place the other on $S^2$.
We can define a minimum separation angle $\theta_{\rm min}$ for kissing points on $S^3$. This 
 angle can be obtained imposing $d_3\geq 1$, using Eq.\  (\ref{d3}):
  \begin{equation}
 \theta_{\rm min} = \left| \theta_j-\theta_i+\Phi_{ij}\right|_{\rm min}=\cos^{-1} \left( \frac{1}{\sqrt{4-d_2^2}}\right)
 \label{theta}
 \end{equation}
 The last expression must be taken with care. When $d_2>\sqrt{3}$ the argument of  $\cos^{-1}$
 is larger than one and makes no formal sense, but tell us that there is not a minimum angle between
 points on different circles separated by that large a distance. Thus, any points on $S^3$ coming from different circles separated a distance $d_3>\sqrt{3}$
 are always kissing points. From the above we easily deduce that if we have a kissing configuration and we add the same constant $c$ to all angles, the resultant 
 configuration remains kissing.

\begin{figure}
\leavevmode
\begin{center}
\includegraphics[angle=0,width=12cm]{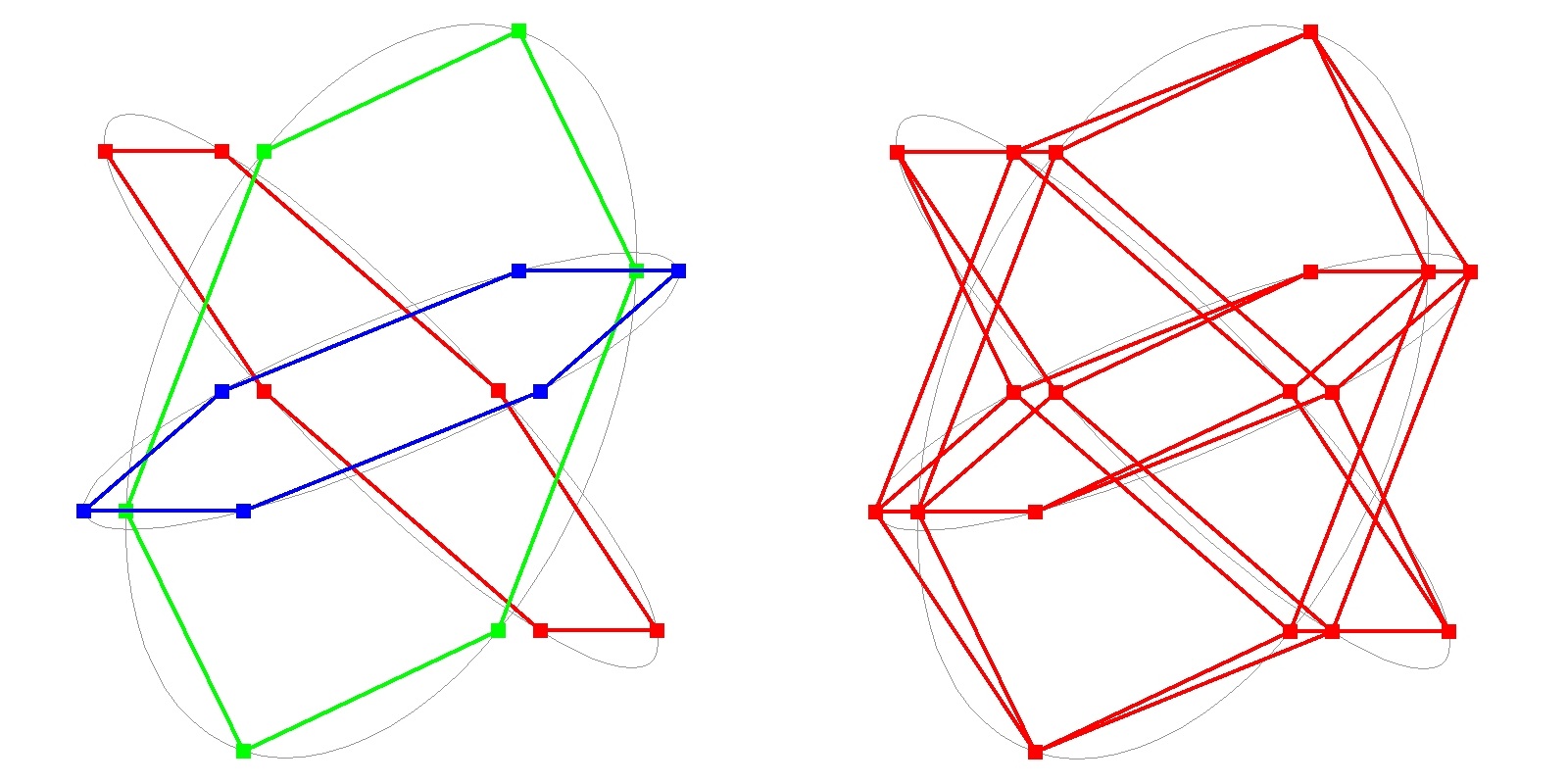}
\caption{3x6 kissing configuration in a orthographic projection to a plane. In the left figure lines join points on the same circle, each circle is drawn in a diferent color. Gray lines 
correspond to full circles. In the right figure lines join nearest neighbors points. After a rotation in 4D, only those points over a line passing through the center are on the same circle, i.e. antipodal points.}
\label{3x6}
\end{center}
\end{figure}

After a rigid body rotation in four dimensions, points on the same circle change to diferent circles. Points only remain over a same circle after a rotation if they are antipodal or, 
in other words, these points have angles separated by $\pi$ radians. Let us have two points on $S^3$ separated by $d_3$ then, their circles on $S^2$  can be separated 
by a maximum distance  $d_{2;max}$ given by:
\begin{equation}
d_{2;max}=\sqrt{4-\left(d_3^2-2\right)^2}.
\end{equation}
 $6\times 4$ kissing configuration (24-cell) after a rotation becomes $12\times 2$ since each circle have two pairs of antipodal points. 
We name a configuration $N\times n$ {\it irreducible} if $N$ it is the maximum integer we can get after any rotation of $S^3$. Any configuration
 with only one point or two antipodal points per circle is irreducible. 
$3\times 6$ kissing configuration is reducible to $9\times 2$ after rotation in $S^3$, see Fig.\ \ref{3x6}.
Anstreicher \cite{An02} showed that the unique antipodal configuration with 24 points is the 24-cell.  Then to see if there are other configurations of 24 kissing points in four dimensions one has to consider configurations of the form $N\times 2$ + $n \times$1 for $N = 11, 10, 9, \ldots , 0$.

\begin{figure}
\leavevmode
\begin{center}
\includegraphics[angle=0,width=12cm]{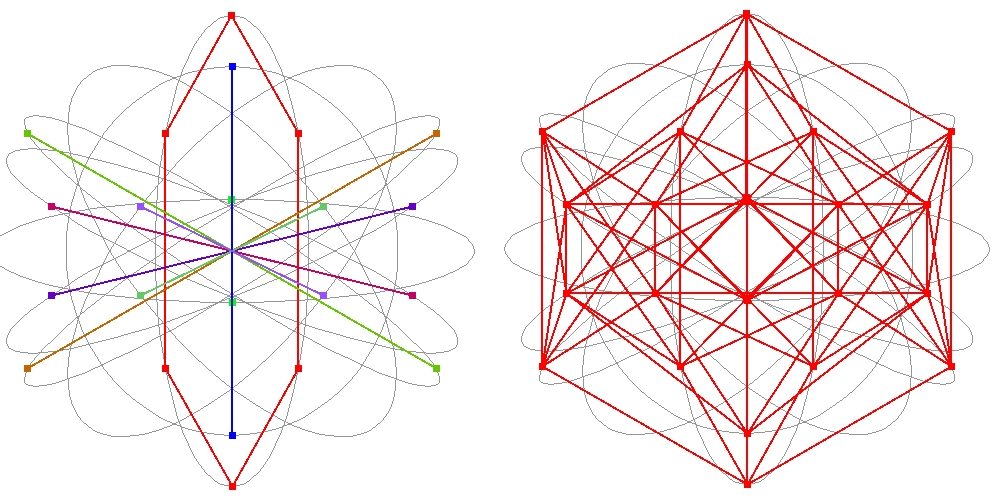}
\caption{Cohn and Woo's 22 spheres kissing configuration in a orthographic projection to a plane. It is shown $1\times 6+6\times 2$ but after a rotation it becomes a $11\times 2$. In the left figure lines join points on the same circle, each circle is drawn in a diferent color. Gray lines 
correspond to full circles. In the right figure lines join nearest neighbors points. }
\label{cohn22}
\end{center}
\end{figure}

Can there exist a 24-point kissing configuration of the form $11\times 2$ + $2\times$ 1?  The only kissing configuration with 11 antipodal points ($11\times 2$)  not simply a subset of the 24-cell known was recently found by Cohn and Woo\cite{CW12}.  This is shown in Figure \ref{cohn22}. Coordinates of this configuration can be obtained from the Hopf map using equation \ref{map3-4} and:

\indent $a=0,\, t=0,$ $\theta=(n-1/2)\frac{\pi}{3}$ $(n=1,2,3,4,5,6)$\\
\indent $a=0\pm i\sqrt{3}/2,$ $t=-1/2,$ $\theta=\pi/2,3\pi/2$\\
\indent $a=\pm\sqrt{3}/2,$ $t=-1/2,$ $\theta=\pi/2,3\pi/2$\\
\indent $a=\pm 2/3 \pm i2/3,$ $t=-1/2,$ $\theta=0,\pi/2$\\

Which corresponds to a $1\times 6+ 8\times 2$ configuration but reducible to a $11\times 2$ after a rotation in 4D.

We now prove that kissing configurations of the form $11\times 2 + 2\times 1$ are not possible.  For say there were such a configuration then by removing each of the singletons as we show below one would get a kissing configuration of the form $12\times 2$--two different kissing configurations with twelve antipodal pairs thus contradicting Anstricher\cite{An02}.

As we mentioned,  each point in $S^3$ can be represented by a point on $S^2$
and an angle. Let say that, for 11$\times$2+2$\times$1, angles are named $\alpha_{i,k}$ where
$i$=1,$\ldots$,13, $k$=1,2 for $i\leq$11 and $k$=1 for $i$=12,13. Antipodal points are
on the same circle and verify: $\alpha_{i,1} = \alpha_{i,2} + \pi$.
So we must check that $\alpha_{12,2} = \alpha_{12,1} + \pi$ is kissing to probe that no irreducible $11\times 2  + 2\times 1$ kissing exists.
Distance in $S^3$ between two points is given by Eq. \ref{d3}, and we can rewrite:
\begin{equation}
d^2_3=2-\sqrt{4-d_{i,j}^2}\cos\left(\alpha_{i,k_i} - \alpha_{j,k_j} + \Phi_{ij}\right),
\end{equation}

where $d_{i,j}$ is the distance on $S^2$ between circles $i$ and $j$ and $\Phi_{i,j}$ is an
angle that depends on relative positions of circles $i$ and $j$ as previously stated.
We remove a point and get $11\times2+1\times1$, since is a kissing config:

\[
2-\sqrt{4-d_{12,j}}\cos\left(\alpha_{12,1} - \alpha_{j,1} + \Phi_{12,j}\right) \geq 1
\]
and
\[
2-\sqrt{4-d_{12,j}}\cos\left(\alpha_{12,1} - \alpha_{j,2} + \Phi_{12,j}\right) \geq 1
\]

for each $j$=1,...,11.

As $\alpha_{j,2} + \pi=\alpha_{j,1}$ we can write the latter expression as:
\[
2-\sqrt{ 4-d_{12,j} }\cos\left(\alpha_{12,1} - \alpha_{j,1}+\pi + \Phi_{12,j}\right) \geq 1
\]
or
\[
2-\sqrt{4-d_{12,j}}\cos\left( (\alpha_{12,1}+\pi) - \alpha_{j,1} + \Phi_{12,j}\right) \geq 1
\]

and easy to get also:
\[
2-\sqrt{4-d_{12,j}}\cos\left( (\alpha_{12,1}+\pi - \alpha_{j,2} + \Phi_{12,j}\right) \geq 1
\]

for each $i$=1,$\ldots$,11, thus, $\alpha_{12,1}+\pi$ (antipodal) is also a kissing point.

\begin{figure}
\leavevmode
\begin{center}
\includegraphics[angle=0,width=10cm]{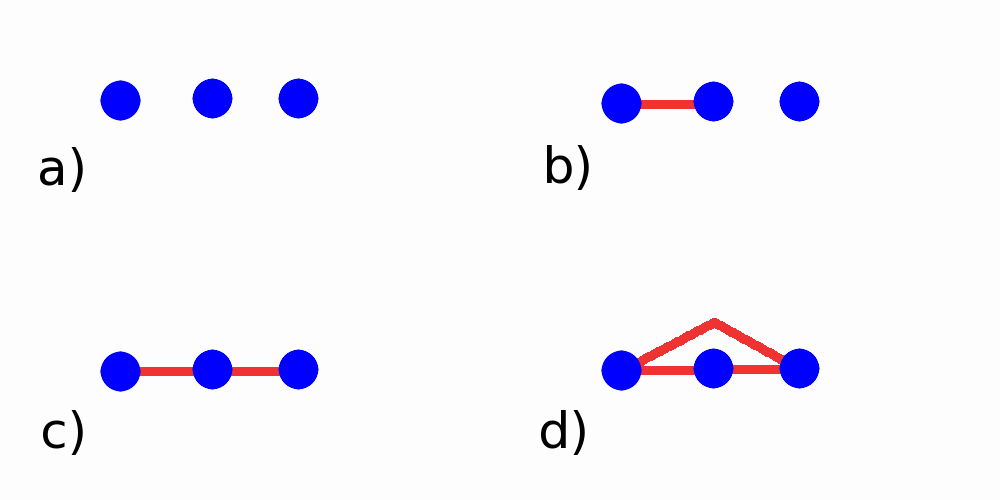}
\caption{All possible covering graphs for 3 spheres. If we have the case a), the we can add 3 antipodal points to that configuration obtaining a not allowed $13\times 2 $ kissing configuration. If we have case b), it is possible to add the antipodal to not bounded sphere getting a $11\times 2+2\times 1$ kissing configuration, which is not allowed. For case c), removing the central sphere, we could add two antipodal points getting a $12\times 2$ kissing config. This must be the 24-cell, thus, not possible. Case d) is not possible because violates the antitransitive property.   }
\label{3diagram}
\end{center}
\end{figure}

Using this proof we also  show that there can be no kissing configurations of the form $10\times 2 + 3\times 1$ 
(and thus certainly no configurations of the form $10\times 2 + 4\times 1$) or of the form
$9\times 2 + 6\times 1$. Indeed if there were a kissing configuration of the form $10\times 2 + 3\times 1$ and an antipodal point were added to one of the singletons one would have a configuration of the form $11\times 2 + 2\times 1$ which we just showed is not possible.  As above if one adds an antipodal point to a $10\times 2 + 3\times 1$ configuration the antipodal point is kissing with the $10\times 2$ antipodal points and of course kissing with its antipodal partner.  The only thing new to be shown is that is kissing with the other two singletons.  
Let us say that a unit sphere $p_1$ is in the cover set of  another unit sphere $p_2$ ($p_1 \in cov(p_2)$) if it is not possible to place a third unit sphere on the antipodal of $p_1$ and get a kissing config. This is clearly symmetric, $p_1$ $\in$ $cov(p_2)$ $\Rightarrow$   $p_2$ $\in$ $cov(p_1)$ and is antitransitive: $p_1$ $\in$ $cov(p_2)$ and $p_2$ $\in$ $cov(p_3)$
$\Rightarrow$   $p_1$ $\notin$ $cov(p_3)$. This can be represented graphically using graphs were triangles are not allowed. Each point in the graph is a unit sphere and a bond between to spheres implies covering. In figure \ref{3diagram} we show the covering posibilities for 3 spheres, what show that a $10\times 2 +3\times 1$ kissing configuration would imply the existence of a forbidden kissing configuration. In these graphs, the existence of a sphere without bonds would imply that we can get a $(N+1)\times 2 +(n-1)\times 1$ kissing configuration from any
$N\times 2 +n\times 1$ kissing arrangement. If we have in the graph an sphere with just one bond, we would be able to obtain a $(N+1)\times 2 +(n-2)\times 1$ from $N\times 2 +n\times 1$.
On the other hand, if we have an sphere $p_1$ with $n/2$ bonds, we can remove all spheres not in $cov(p_1)$ and then after adding antipodal to spheres in $cov(p_1)$ we get a $(N+n/2)\times 2$ kissing configuration from $N\times 2+ n\times 1$.

\begin{figure}
\leavevmode
\begin{center}
\includegraphics[angle=0,width=6cm]{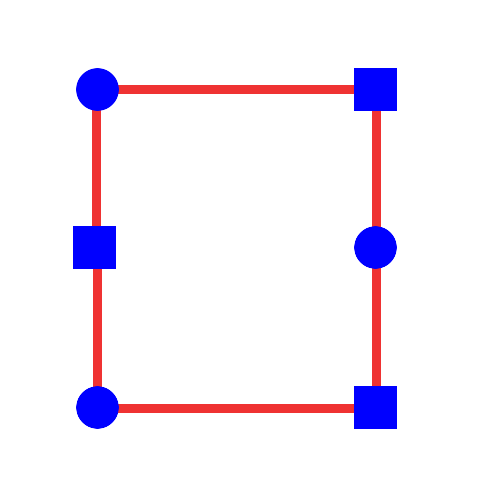}
\caption{Covering graph for 6 spheres. Removing squares and adding antipodal to circles we get a $12\times 2$ kissing configuration or removing circles and adding antipodal to squares we get a different $12\times 2$ kissing configuration.}
\label{6diagram}
\end{center}
\end{figure}

Next let us show that kissing configurations of the form $9\times 2 + 6 \times 1$ cannot exist.  
If we have any of the six spheres in the singleton with a number of elements in its cover set less or greater than 2, then we would get a forbidden kissing config as demostrated above.
Then, the only graph to analyze is shown in Fig. \ref{6diagram}, that would imply the existence of two diferent $12\times2$ kissing configurations.
We have not been able  to continue this line of potential proof of the uniqueness of the 24-cell to  configurations of the form $8\times 2 + 8\times 1$.

But we are able to delimit the possible maximal kissing configurations in four dimensions to deriving from at least sixteen circles on $S^2$.  In trying to find a $16\times 1$ configuration it is easy to construct one analytically starting from $3\times 5+1\times 1$ where
first 3 circles are equispaced on the equator and the last one is on
the north pole. We place, for example, $\theta$=0,61,122,185,250 degrees for
circles on equator and $\theta$=300 for the circle at the pole and
that config is kissing and after a rotation in 4D becomes a $16\times 1$ since
no point is antipodal of any other. If we put $\theta =n\pi/3 $, $n=0,1,2,3,4$
the config is also kissing but after a rotation it changes to an
irreducible $6\times 2+4\times 1$.
The configuration of the form n x 1 (n the number of circles on $S^2$) with the largest n of which we are aware is
has $N$= 22\cite{Sl94}.
While we have not been able to prove that the 24-cell is the unique kissing configuration in four dimensions, for the first time we have been able to delimit the space of other configurations that could possibly be kissing. We hope that our findings and approach may be helpful in learning more about kissing configurations in four and higher dimensions.

\begin{thebibliography} {99}

\bibitem{Me98}  Melisen JBM . How Different Can Colours Be? Maximum Separation of Points on a Spherical Octant.  {\it 
Proceedings: Mathematical, Physical and Engineering Sciences} {\bf 454} (1973): 1499–1508. (1998).

\bibitem{CS98} J. H. Conway and N. J. A. Sloane, Sphere Packing, Lattices and Groups  Springer-Verlag NY 3rd edition (1998).

\bibitem{Ta30}Tammes, R.M.L. On the Origin Number and Arrangement of the Places of Exits on the
 Surface of Pollengrains,{\it Rec. Trv. Bot. Neerl.} {\bf 27,}, 1--84 (1930).

\bibitem{BG03}[Bruinsma RF, Gelbart WM, Reguera D, Rudnick J, Zandi R. "Viral Self-Assembly as a Thermodynamic Process". {\it Physical Review Letters} {\bf 90} (24): 248101–1–248101–4. (1990).

\bibitem{LD03} T. Liu, E. Diemann, H. Li, A. W. Dress, and A. Muller, A. {\it Nature}  {\bf  426}, 59 (2003).

\bibitem{MR08} J. Mikhael,  J. Roth,  L. Helden and  C. Bechinger, Archimedean-like tiling on decagonal quasicrystalline surfaces  {\it Nature} {\bf 454}, 501-504 (2008).

\bibitem{BM12}T. Bohlein, J. Mikhael and C. Bechinger,
Observation of kinks and antikinks in colloidal monolayers driven across ordered surfaces 
{\it Nature Materials} {\bf 11}, 126–130 (2012).

\bibitem{Ca11}  Casselman, W. The Difficulties of Kissing in Three Dimensions.
{\it Notices of the AMS} {\bf 51,} 884 (2004).

\bibitem{SW53}   Sch\"utte, K.  \&    van der Waerden, B. L.
Das Problem der dreizehn Kugeln. 
{\it Math.\ Ann.\/} {\bf 125,} 325--334  (1953).

\bibitem{De72} Delsarte, P. Bounds for unrestricted codes, by linear
programming, {\it Philips Res. Rep.} {\bf 27,} 272--289  (1972).

\bibitem{De73}  Delsarte, P. An algebraic approach to the association
schemes of coding theory, {\it Philips Res. Rep. Suppl.} vi+97 (1973).

\bibitem{DG77}   Delsarte, P. Goethals, J. M. \& Seidel, J.J. Spherical
codes and designs, {\it  Geom. Dedicata} {\bf 6,} 
363--388  (1977).

\bibitem{Le79} Levenshtein, V. I. 
On bounds for packing in $n$-dimensional Euclidean space.
{\it Sov.\ Math.\ Dokl.\/} {\bf 20(2),} 417--421 (1979).  


\bibitem{OS79}  Odlyzko, A. M.   \&     Sloane,  N. J. A.
New bounds on the number of unit spheres that can touch a unit sphere in $n$ dimensions.
{\it J. of Combinatorial Theory } {\bf A26,} 210--214 (1979). 

\bibitem{AB00} Arestov, V. V.  \& Babenko, A. G.  Estimates for the maximal
value of the angular code distance for 24 and
25 points on the unit sphere in R4, {\it Math. Notes} {\bf 68}
419--435 (2000).

\bibitem{Mu03} Musin,   O. R.
The problem of the twenty-five spheres. Russ. Math. Surv. {\bf 58,} 794--795 (2003). 

\bibitem{BV09}  Bachoc, C.  \&  Vallentin, F.
New upper bounds for kissing numbers from semidefinite programming, {
\it J. Amer. Math. Soc.} {\bf 21,} 909--924 (2008).

\bibitem{MV09} Mittelmann, H. D. \& Vallentin, F. 
High accuracy semidefinite programming bounds for kissing numbers.
Experiment. Math. {\bf 19,} 174-178 (2010).



\bibitem{Ho31} Hopf, H.
{\it \"Uer die Abbildung der Dreidimensionalen Sphare auf die
 Kugelfl\"ahe. }
  {\it Math.\ Ann.} {\bf 104}  (1931)
  Reprinted in {\it Selecta Heinz Hopf },   p. 38--63,
  Springer-Verlag Berlin Heidelberg New York, (1964).
  
  
\bibitem{AP07}Altschuler, E. L. \& \ P\'erez-Garrido, A.
Symmetric four-dimensional polytope and visualization method in four, eight and sixteen dimensions using Hopf maps. 
Phys.\ Rev.\ E {\bf 76,} 016705-1 (2007).




 \bibitem{An02} Anstreicher, K. M. Improved Linear Programming Bounds for Antipodal Spherical Codes, {\it Dicrete and Computational Geometry}
 { \bf  8,} 107--114 (2002).

\bibitem{CW12} Cohn, H. \& Woo, H. Three-point bounds for energy minimization. {\it Journal of the American Mathematical Society} {\bf 25,} 929-958 (2012)


\bibitem{Sl94} Sloane, N. J. A. {\it Tables of Spherical Codes.}
\verb! http://www2.research.att.com/~njas/packings/!



\end{thebibliography}
\end{document}